\documentclass[a4paper,10pt,twoside]{article}
\usepackage[TS1,T1]{fontenc}
\usepackage{enumitem}
\usepackage[toc]{appendix}
\usepackage{float}
\usepackage{textcomp}
\usepackage{amsmath,amsfonts,verbatim,afterpage,euscript,mathrsfs,amssymb}
\usepackage{amsthm}
\usepackage{dsfont}
\usepackage{hyperref}
\usepackage{pst-fill,pst-grad}
\usepackage{textcomp}
\usepackage{multicol}
\newcommand{\mysection}{\setcounter{equation}{0} \section}

\newtheorem{Theorem}{Theorem}
\newtheorem{Lemma}{Lemma}
\newtheorem{Proposition}{Proposition}

\newtheorem{Remark}{Remark}

\def \vu{\vec{u}}

\def \R{\mathbb{R}}
\def \M{\mathcal{M}}
\def \m{\mathfrak{m}}
\title{\bf  Global weak solutions for a variation of the Whitham equation}
\usepackage{authblk}
\author{Diego Chamorro\footnote{\emph{diego.chamorro@univ-evry.fr} }}
\affil{\footnotesize LaMME, Univ. Evry, CNRS, Universit\'e Paris-Saclay, 91025, Evry, France.}
\author{Mar\'ia Eugenia Mart\'inez\footnote{\emph{maria.martinez.m@uchile.cl} } \thanks{M.E.M. was partially funded by Chilean research grant ANID Exploraci\'on 13220060, has been partially supported by the project CRISIS (ANR-20-CE40-0020-01), operated by the French National Research Agency (ANR), and by the HiCE project, operated by the Universit\'e Claude Bernard Lyon 1. This work was supported by Centro de Modelamiento Matem\'atico (CMM) BASAL fund FB210005 for center of excellence from ANID-Chile.}}
\affil{\footnotesize Universidad de Chile, Chile \& Universit\'e Claude Bernard Lyon 1, France.}

\begin{document}
\sloppy
\maketitle
\begin{scriptsize}
\abstract{We study in this article a variation of the Whitham equation which was introduced as an alternative to the KdV equation. We first prove the global existence of weak solutions, then we establish a regularity criterion from which we deduce the uniqueness of weak solutions. Local in time criterions for regularity and uniqueness are also given.}\\[3mm]
{\bf \scriptsize Keywords: Whitham equation; global weak solutions; regularity criterion; uniqueness.}\\
\textbf{\scriptsize Mathematics Subject Classification: 35Q35; 35R11.}
\end{scriptsize}
\mysection{Introduction}

In this article  our goal is to study the following \emph{modified} Whitham equation:
\begin{equation}\label{EqIntro}
\begin{cases}
\partial_t u = -(-\Delta)^{\frac{1}{2}}\M(u)+\partial_x\left(\dfrac{u^2}{2}\right),\\[3mm]
u(0,x)=u_0(x), \quad x\in \R,
\end{cases}
\end{equation}
where $u_0 : \R \longrightarrow \R$ is an initial data such that $u_0 \in L^2(\R)$ and the operator $\M$ is defined in the Fourier-level by the following expression
\begin{equation}\label{Def_M}
\widehat{\M (u)}(\xi):= \mathfrak m(\xi)\widehat u (\xi), \quad \text{and} \quad 
\mathfrak m(\xi):=
\sqrt{(1+\xi^2)\frac{\tanh(\xi)}{\xi}}.
\end{equation}
Note that the symbol of the operator $\M$, given by the positive function $\m$, satisfies $\m(0)=1$ and $\m(\xi)\sim |\xi|^{\frac{1}{2}}$ if $|\xi|\gg 1$. Recall that the action of the fractional power of the Laplacian $(-\Delta)^{\frac{1}{2}}$ is given in the Fourier variable by $\widehat{(-\Delta)^{\frac{1}{2}}f}=|\xi|\widehat{f}(\xi)$. Note also that the operators $(-\Delta)^{\frac{1}{2}}$ and $\mathcal{M}$ commute (as it can be easily checked in the Fourier variable) and we have $(-\Delta)^{\frac{1}{2}}\M(u)=\M((-\Delta)^{\frac{1}{2}}u)$.\\

Let us point out now that the equation (\ref{EqIntro}) is a variation of the capillary-gravity Whitham equation, given by 
\begin{equation}\label{EqIntroWhitham}
\begin{cases}
\partial_t u=-\partial_x\left(\M (u) - \dfrac{u^2}2\right),\\[3mm]
u(0,x)=u_0(x), \quad x\in \R,
\end{cases}
\end{equation}
This problem was first introduced by Whitham \cite{whitham} as a shallow water model that, unlike KdV, would allow phenomena like wave breaking and peaking. The Whitham equation without surface tension, proposed by Whitham in \cite{whitham}, reads as follows 
\begin{equation}\label{GeneralWhithamEquation}
    \partial_t u + \partial_x \left( \mathcal K_{h_0} u + \frac34\frac{c_0}{h_0} u^2\right)=0,  \quad 
    \widehat{\mathcal K_{h_0}u}(\xi):= \sqrt{\frac{g\tanh(h_0\xi)}{\xi}},
\end{equation}
where $c_0 = \sqrt {gh_0}$ is the limiting long-wave speed, $h_0$ is the undisturbed fluid depth and $g$ is gravitational acceleration. When considering the KdV model
\[
\partial_t u + \partial_x \left( c_0 u + \frac16 c_0 h_0^2\partial_x^2 u + \frac34 \frac{c_0}{h_0} u^2\right)=0,
\]
one finds that the phase velocity of the KdV equation, given by $c_{\text{KdV}}(\xi)=c_0 - \frac16 c_0h_0^2 \xi^2$, is a second-order approximation to the phase velocity of the Whitham equation \eqref{GeneralWhithamEquation}:
\[
c_\text{W}(\xi)= \sqrt{\frac{g\tanh(h_0\xi)}{\xi}},
\]
in the low-frequency (or long-wave) limit, where $\vert \xi h_0\vert\ll 1$. However, for high frequencies, the quantity $c_{\text{KdV}}$ fails to approximate the ``full dispersion relation'' $c_{\text W}$ effectively. This discrepancy is tied to the existence of solitary and periodic waves in the KdV model, a phenomenon absent in traditional shallow water theory. As Whitham explained in \cite{whitham2}, in his effort to create a model that admits wave peaking and breaking, he proposed the \emph{full dispersion equation}, or Whitham equation with no surface tension \eqref{GeneralWhithamEquation}. This model successfully combines the full dispersion relation $c_{\text W}$, which is more moderate than KdV, with the long-wave nonlinearity of classical shallow water theory.\\

As subsequent studies showed, Whitham's intuition about wave breaking (where solutions remain bounded but their derivatives become unbounded in finite time) and peaking was at least partially correct. It was later shown that the Whitham equation without surface tension  \eqref{GeneralWhithamEquation} exhibits wave breaking (see  \cite{Hur} and \cite{KleinSaut} for numerical analysis), although this phenomenon may not necessarily reflect the behavior of actual water waves (see \cite{KleinLinaresPilodSaut}). Equation \eqref{GeneralWhithamEquation} also admits the existence of periodic waves with sharp crests (as shown in \cite{EhrnstromWahlen}), periodic traveling waves \cite{EhrnstromKalisch} and solitary waves \cite{EhrnstromGrovesWahlen, SW}. Of course, the dynamics change when surface tension is considered: indeed, for \eqref{EqIntroWhitham}, the high-frequency behavior of the dispersion becomes stronger, resembling that of the KdV equation and this results in different dynamics and more ``well-behaved'' solutions. Finally, we refer to \cite{Borluka}, where the authors study numerically the viability of the equation proposed by Whitham \eqref{GeneralWhithamEquation} as a water waves model. \\

Regarding the well-posedness of the Whitham equation, it has been established that it is locally well-posed in the space $H^s(\R)$ for $s > 3/2$ (see \cite{EhrnstromEscherPei, KleinLinaresPilodSaut}). Local (and global) existence for a class of dispersive perturbations to the Burgers equation was studied in \cite{MolinetPilodVento}. However, to date, there are no known global well-posedness results for \eqref{EqIntroWhitham}, particularly in the natural energy space, which corresponds to $H^{\frac{1}{4}}$, as we will see below. \\

In this article, we will first address the problem of \emph{global existence} and we shall study a modified Whitham equation, where the derivation $\partial_x$ is replaced by the square root of the Laplacian as in equation (\ref{EqIntro}) above. Let us remark that this modification will introduce an interesting (weak) diffusion effect whose action will be compensated by the nonlinear term $\partial_x\left(\dfrac{u^2}{2}\right)$. Indeed, since $\mathfrak{m}(\xi)\simeq 1$ if $|\xi|<1$, then the action of the operator $(-\Delta)^{\frac{1}{2}}\M$ in the Fourier level is essentially given by the Fourier multiplier $|\xi|$ whose regularizing effect in this regime of small frequencies is in competition with the derivative present in the nonlinearity: thus the regularizing effect of the term $(-\Delta)^{\frac{1}{2}}\M$ is not necessarily strong enough to produce in this system a clear and sharp gain of regularity.\\

In our first result we construct global weak solutions for the equation (\ref{EqIntro}).
\begin{Theorem}[\bf Global Solutions]\label{Theorem_GlobalExistence}
For all initial data $u_0\in L^2(\R)$ the equation (\ref{EqIntro}) admits a global weak solution $u\in L^\infty_t(L^2_x)$ that satisfies the following energy inequality
$$\|u(t,\cdot)\|^2_{L^2}+\int_{0}^{+\infty}\|u(s,\cdot)\|_{\mathcal{N}}^2ds\leq \|u_0\|_{L^2}^2,$$
where the functional $\|\cdot\|_{\mathcal{N}}$ is defined in the expression (\ref{Def_FunctionalN}) below.
\end{Theorem}
Let us remark that the presence of the operator $(-\Delta)^{\frac{1}{2}}\M$ in the equation (\ref{EqIntro}) will naturally lead us (when studying energy inequalities, see Proposition \ref{Propo_Energy} below) to the functional space $\mathcal{N}(\R)$ defined as the set of functions $f:\R\longrightarrow \R$ such that 
\begin{equation}\label{Def_FunctionalN}
\|f\|_{\mathcal{N}}=\left(\int_{\R}\left|(-\Delta)^{\frac{1}{4}}\M^{\frac{1}{2}}(f) \right|^2dx\right)^{\frac{1}{2}}<+\infty,
\end{equation}
where the symbol in the Fourier level of the operator $(-\Delta)^{\frac{1}{4}}\M^{\frac{1}{2}}$ is given by $|\xi|^\frac{1}{2}\m^{\frac12}(\xi)$. Note that, by the Plancherel theorem, we have $f\in \mathcal{N}(\R)$ if 
$$\left(\int_{\R}|\xi|\left((1+\xi^2)\frac{\tanh(\xi)}{\xi}\right)^{\frac{1}{2}}|\widehat{f}(\xi)|^2d\xi\right)^{\frac{1}{2}}<+\infty,$$
and as we shall see in the pages below this functional space $\mathcal{N}(\R)$ conveys some regularity information that can be stated in terms of the usual Sobolev spaces (see Lemma \ref{Lemm_InterpolMain}). \\

However, in spite of the presence of the operator $(-\Delta)^{\frac{1}{2}}\M$ and of the Sobolev-like information given by the space $\mathcal{N}(\R)$ which will be induced by this operator, the nonlinear term $\partial_x(u^2)$ must be treated carefully when considering existence results and to establish the previous theorem we will proceed in several steps. Indeed, first we will introduce a hyperviscosity perturbation, that will provide us with some extra regularity and then we will perform a fixed point argument in the space $L^\infty_t(L^2_x)$. Next we will prove that the solutions of this perturbed system are regular and with this information we will establish an energy inequality that will allow us to deduce uniform estimates. This energy inequality will be helpful to extend the time of existence of the solutions and to get rid of the hyperviscosity perturbation in order to recover the original system. Let us remark here that, although we obtain global weak solutions to the system (\ref{EqIntro}), in our proof we need to consider a limit process (up to subsequences) and therefore we cannot grant the uniqueness of such solutions.\\

Let us remark that most of the previous steps can be applied to the equation (\ref{EqIntroWhitham}), such as the hyperviscosity perturbation and the global existence results for the mollified problem. However, without a suitable \emph{a priori} estimate, we can not recover a global weak solution for the original system (\ref{EqIntroWhitham}).\\

Another problem related to this type of solutions is their regularity. This problem is of course interesting by itself, but as we shall see later on, it can help to solve the uniqueness issue mentioned above.  As pointed out previously, the regularity properties of the operator $(-\Delta)^{\frac{1}{2}}\M$ (which behave in the Fourier level as $|\xi|^{\frac{3}{2}}$ if $|\xi|\gg 1$ but only as $|\xi|$ if $|\xi|<1$) seem too weak to compensate the effect of the nonlinear drift term $\partial_x(u^2)$ and thus some additional hypotheses will be needed to obtain some regularity results.\\

 In our next result we give a regularity criterion: 
\begin{Theorem}[\bf Regularity Criterion]\label{Theorem_SerrinLinfty}
Consider a initial data $u_0\in L^2(\R)$ such that $u_0\in \dot{H}^{\rho}(\R)$ for some $\rho>0$ large enough and consider an associated weak solution $u\in L^\infty_t(L^2_x)\cap L^2_t(\mathcal{N}_x)$ to the equation (\ref{EqIntro}). Assume moreover that $u\in L^\infty_t(L^\infty_x)$, then the solution $u$ belongs to the space $L^\infty_t(\dot{H}^{\rho}_x)$. 
\end{Theorem}
\noindent Of course, once we have enough regularity, we can study the uniqueness for the weak solutions obtained in Theorem \ref{Theorem_GlobalExistence} and in this sense we have our last result:
\begin{Theorem}[\bf Uniqueness]\label{Theorem_Uniqueness}
Consider a initial data $u_0\in L^2(\R)$ such that $u_0\in \dot{H}^{\rho}(\R)$ for some $\rho>0$ large enough and consider an associated weak solution $u\in L^\infty_t(L^2_x)\cap L^2_t(\mathcal{N}_x)$ to the equation (\ref{EqIntro}). Assume moreover that $u\in L^\infty_t(L^\infty_x)$, then the solution $u$ is unique.
\end{Theorem}
To finish this introduction, let us make some remarks. From the point of view of the competition between the regularizing effects of the operator $(-\Delta)^{\frac{1}{2}}\M$ and the nonlinearity $\partial_x(u^2)$, the equation (\ref{EqIntro}) is almost critical in the sense that the operator $\M$ has very mild regularizing properties. This scenario can be compared to the study of the quasi-geostrophic  equation, which has some additional structure but is very difficult to study and still contains open problems, see \emph{e.g.} \cite{Dab} and the references therein for more details. Let us observe also that the regularity criterion given above can be obtained due to the properties of the semi-group associated to the operator $(-\Delta)^{\frac{1}{2}}\M$. If we consider the operator $\partial_x\M$, as in the system (\ref{EqIntroWhitham}), then the study of the regularity of global solutions would probably require a completely different treatment. Finally, let us note that the regularity criterion given in Theorem \ref{Theorem_SerrinLinfty} can be relaxed when considering a local in time framework.\\

The plan of the article is the following. In Section \ref{Secc_Notation}, we introduce some of the notation as well as some properties that will be useful here. In Section \ref{Secc_GlobalWeak}, we will prove Theorem \ref{Theorem_GlobalExistence} and in Section \ref{Secc_Regularity}, we study the Theorem \ref{Theorem_SerrinLinfty}. In Section \ref{Secc_Uniqueness}, we apply the previous results to the uniqueness problem and we establish Theorem \ref{Theorem_Uniqueness}. Finally, in Section \ref{Secc_Local}, we relax (in a bounded in time setting) the boundedness hypothesis used in Theorem \ref{Theorem_SerrinLinfty} and  \ref{Theorem_Uniqueness}.
\mysection{Notations and useful estimates}\label{Secc_Notation}
We will denote by  $L^\infty_t(L^2_x)$ the Lebesgue space in time and space given by the condition
$$\|u\|_{L^\infty_t(L^2_x)}=\underset{0<t<T}{\sup}\left(\int_{\R}|u(t,x)|^2dx\right)^{\frac{1}{2}}=\underset{0<t<T}{\sup}\|u(t,\cdot)\|_{L^2}<+\infty,$$
where $u:[0,+\infty[\times \R\longrightarrow \R$.\\

\noindent For $\sigma>0$ the homogeneous Sobolev spaces will be denoted by $\dot{H}^{\sigma}(\R)$ and for a regular function $f:\R\longrightarrow \R$ we have
$$\|f\|_{\dot{H}^{\sigma}}=\|(-\Delta)^{\frac{\sigma}{2}}f\|_{L^2}.$$
Recall that we have the Sobolev inequalities
$$\|f\|_{L^q}\leq C\|f\|_{\dot{H}^{\sigma}}, \quad \mbox{with } 0<\sigma<\frac{1}{2}  \quad \mbox{and } \frac{1}{q}=\frac{1}{2}-\sigma.$$
Moreover, we also have the Sobolev embedding $L^2(\R)\cap \dot{H}^{\sigma}(\R)\subset L^\infty(\R)$ if $\sigma>\frac{1}{2}$. See \cite[Theorem 6.2.4]{Grafakos} for a proof of these inequalities.\\

When dealing with Sobolev spaces and a product of two functions, we have at our disposal two useful results. Indeed if $f,g:\R\longrightarrow \R$ are two functions such that $f,g\in \dot{H}^{\sigma}(\R)$ with $\sigma>0$ and $f,g\in \dot{H}^{\delta}(\R)$ with $0\leq \delta\leq \frac{1}{2}$, then we have the following law product
\begin{equation}\label{LawProduct_Ineq}
\|fg\|_{\dot{H}^{\sigma+\delta-\frac{1}{2}}}\leq C\left(\|f\|_{\dot{H}^{\sigma}}\|g\|_{\dot{H}^{\delta}}+\|g\|_{\dot{H}^{\sigma}}\|f\|_{\dot{H}^{\delta}}\right),
\end{equation}
see \cite[Theorem 4.1]{PGLR1} for a general proof of this fact. Note that we also have the following  Kato-Ponce inequality (also know as the fractional Leibniz rule):
\begin{equation}\label{KatoPonce_Ineq}
\|(-\Delta)^{\frac{\sigma}{2}} fg\|_{L^2}\leq C\left(\|(-\Delta)^{\frac{\sigma}{2}} f\|_{L^2}\|g\|_{L^\infty}+\|(-\Delta)^{\frac{\sigma}{2}} g\|_{L^2}\|f\|_{L^\infty}\right),
\end{equation}
see \cite{GrafakosOh} for a proof of this result in a more general framework.

\mysection{Global weak solutions}\label{Secc_GlobalWeak}
In this section we will obtain global in time solutions for the system (\ref{EqIntro}). We will first introduce a hyperviscosity modification of the equation and then we will obtain weak solutions by a suitable limiting process.\\

First, let us fix a parameter $\varepsilon>0$ and consider the operator $\mathcal{L}=(-\Delta)(Id-\Delta)$ whose action in the Fourier level reads as follows
\begin{equation}
\widehat{\mathcal{L}(u)}(\xi)=|\xi|^2\left(1+|\xi|^2\right)\widehat{u}(\xi).
\end{equation}
Next, we consider a mollifying function $\varphi_\varepsilon(x)=\frac{1}{\varepsilon}\varphi(\frac{x}{\varepsilon})$, where $\varphi\geq 0$ is a test function such that $\displaystyle{\int_\mathbb{R}\varphi(x)dx=1}$. We thus obtain 
\begin{equation}\label{EqMollified}
\begin{cases}
\partial_t u = -\varepsilon \mathcal{L}(u)-(-\Delta)^{\frac{1}{2}}\M(u)+\partial_x\left(\frac{u^2}{2}\right),\\[3mm]
u(0,x)=\varphi_\varepsilon\ast u_0(x), \quad x\in \R,
\end{cases}
\end{equation}
Note that, at least formally, if we set $\varepsilon \to 0$ we recover from the previous equation the original system (\ref{EqIntro}).\\

Remark now that for all $\varepsilon>0$, the operator $-\varepsilon \mathcal{L}$ generates a convolution semi-group given by $e^{-\varepsilon t\mathcal{L}}$ with $t>0$ (\emph{i.e.} we have $e^{-\varepsilon (t+s)\mathcal{L}}=e^{-\varepsilon t\mathcal{L}}e^{-\varepsilon s\mathcal{L}}$). In particular, for a function $\varphi\in \mathcal{S}(\mathbb{R})$ we have 
\begin{equation}\label{Semigroupe}
\widehat{e^{-\varepsilon t\mathcal{L}}\varphi}(\xi)=e^{-\varepsilon  t(|\xi|^{2}(1+|\xi|^2))}\widehat{\varphi}(\xi),
\end{equation}
which allows to write $e^{-\varepsilon t\mathcal{L}}\varphi(x)=\mathcal{G}_{\varepsilon t}\ast \varphi(x)$, where the kernel $\mathcal{G}_{\varepsilon t}$ is such that $\widehat{\mathcal{G}_{\varepsilon t}}(\xi)=e^{-\varepsilon  t(|\xi|^{2}(1+|\xi|^2))}$. Note that we have $\|\mathcal{G}_{\varepsilon t}\|_{L^1}\leq C$ and since we have 
$$\widehat{e^{-\varepsilon t\mathcal{L}}\varphi}(\xi)=e^{-\varepsilon  t(|\xi|^{2}(1+|\xi|^2))}\widehat{\varphi}(\xi)=e^{-\varepsilon  t|\xi|^{2}}e^{-\varepsilon  t|\xi|^4}\widehat{\varphi}(\xi),$$
we can write 
$$e^{-\varepsilon t\mathcal{L}}\varphi=\mathcal{G}_{\varepsilon t}\ast \varphi=\mathfrak{g}_{\varepsilon t}\ast \mathfrak{h}_{\varepsilon t}\ast \varphi,$$
where $\mathfrak{g}_{\varepsilon t}$ is the usual heat kernel and $\mathfrak{h}_{\varepsilon t}$ is the kernel associated to the semi-group $e^{-\varepsilon  t(-\Delta)^2}$. Of course we have $\|\mathfrak{h}_{\varepsilon t}\|_{L^1}\leq C$ and by a straightforward computation in the Fourier level we also have $\|\partial_x\mathfrak{h}_{\varepsilon t}\|_{L^2}\leq C(\varepsilon t)^{-\frac{3}{8}}$.\\

With these preliminaries, using the Duhamel formula we can consider the following integral representation formula of (\ref{EqMollified}):
\begin{equation}\label{Duhamel}
u(t,x)=\underbrace{\mathcal{G}_{\varepsilon t}\ast (\varphi_\varepsilon\ast u_0)(x)}_{(1)}-\underbrace{\int_{0}^{t}\mathcal{G}_{\varepsilon (t-s)}\ast(-\Delta)^{\frac{1}{2}}\M(u)ds}_{(2)}+\underbrace{\int_{0}^{t}\mathcal{G}_{\varepsilon (t-s)}\ast\partial_x\left(\frac{u^2}{2}\right)ds}_{(3)},
\end{equation}
and we will apply a fixed point argument to the previous equation (\ref{Duhamel}) in the space $L^\infty([0,T], L^2(\R))$ for some time $T>0$ to be fixed later. 
\begin{Proposition}\label{Propo_ExistenceApprox}
Consider a initial data such that $u_0\in L^2(\R)$. For a fixed parameter $\varepsilon>0$, the equation (\ref{Duhamel}) admits a mild solution $u_\varepsilon\in L^\infty([0,T], L^2(\R))$ where the time $T$ depends on the initial data $u_0$ and on the parameter $\varepsilon$.
\end{Proposition}
\noindent {\bf Proof.} As announced, we will perform a fixed point argument to obtain a solution of the equation (\ref{Duhamel}). Indeed, we have the following result: 
\begin{Lemma}\label{Lemma_Existence}
Let $(E, \|\cdot\|_E)$ be a Banach space. Consider $e_0\in E$ such that  $\|e_0\|_E\leq \delta$ and consider $L:E\longrightarrow E$ and $B:E\times E\longrightarrow E$ two bounded linear and bilinear applications, \emph{i.e.}:
$$\|L(e)\|_E\leq C_L\|e\|_E\qquad \mbox{and}\qquad \|B(e,f)\|_E\leq C_B\|e\|_E\|f\|_E,$$
for all $e,f\in E$. If we have the following relationship between the constants $\delta, C_L$ and $C_B$
\begin{equation}\label{ConstantesContinuiteBeirao}
0<3C_L<1,\qquad 0<9C_B\delta <1\qquad\mbox{and} \qquad C_L+6C_B\delta <1,
\end{equation}
then the equation
$$e=e_{0}+L(e)-B(e,e),$$
admits a unique solution $e\in E$ such that $\|e\|_E\leq 3\delta$.
\end{Lemma}
\noindent See a proof of this lemma in \cite{ChNS}. Thus, in order to apply this result, we need to estimate the terms (1), (2) and (3) given in (\ref{Duhamel}) above in the functional space $L^\infty_t(L^2_x)$. 
\begin{itemize}
\item For the initial data, by the Minkowski inequalities for the convolution, we easily obtain
\begin{equation}\label{DatoInicial}
\|\mathcal{G}_{\varepsilon t}\ast (\varphi_\varepsilon\ast u_0)\|_{L^2}\leq \|\mathcal{G}_{\varepsilon t}\|_{L^1}\|\varphi_\varepsilon\|_{L^1}\|u_0\|_{L^2}\leq C\|u_0\|_{L^2}.
\end{equation}
\item For the linear term (2), by duality and by the properties of the operators $(-\Delta)^{\frac{1}{2}}$  and $\mathcal{M}$ we write
\begin{eqnarray*}
\left\|\int_{0}^{t}\mathcal{G}_{\varepsilon (t-s)}\ast(-\Delta)^{\frac{1}{2}}\M(u)ds\right\|_{L^2}&=&\underset{\|\psi\|_{L^2}\leq 1}{\sup}\left|\int_{\R}\left(\int_{0}^{t}\mathcal{G}_{\varepsilon (t-s)}\ast(-\Delta)^{\frac{1}{2}}\M(u)ds\right)\psi dx\right|\\
&=&\underset{\|\psi\|_{L^2}\leq 1}{\sup}\left|\int_{0}^{t}\int_{\R} [(-\Delta)^{\frac{1}{2}}\M(\mathcal{G}_{\varepsilon (t-s)})\ast u]\psi dxds\right|,
\end{eqnarray*}
from which we deduce, by the properties of the convolution and by the Cauchy-Schwarz inequality in time and space 
\begin{eqnarray*}
\left\|\int_{0}^{t}\mathcal{G}_{\varepsilon (t-s)}\ast(-\Delta)^{\frac{1}{2}}\M(u)ds\right\|_{L^2}\leq & \underset{\|\psi\|_{L^2}\leq 1}{\sup}\|u\|_{L^2_t(L^2_x)}\left\|[(-\Delta)^{\frac{1}{2}}\M(\mathcal{G}_{\varepsilon (t-s)})] \ast\psi\right\|_{L^2_t(L^2_x)}.
\end{eqnarray*}
At this point, we remark that
\begin{equation}\label{EstimationNoyau2}
\left\|[(-\Delta)^{\frac{1}{2}}\M(\mathcal{G}_{\varepsilon (t-s)})] \ast\psi\right\|_{L^2_t(L^2_x)}\leq \frac{C}{\sqrt{\varepsilon}}\|\psi\|_{L^2}.
\end{equation}
Indeed, by the Plancherel formula and by the Fubini theorem we write
$$\left\|[(-\Delta)^{\frac{1}{2}}\M(\mathcal{G}_{\varepsilon (t-s)})] \ast\psi\right\|_{L^2_t(L^2_x)}^2\leq \int_{0}^{t}\int_{\R}|\xi|^2\left((1+\xi^2)\frac{\tanh(\xi)}{\xi}\right)e^{-2\varepsilon(t-s)|\xi|^2\left(1+|\xi|^2\right)}|\widehat{\psi}(\xi)|^2d\xi ds.$$
Thus by a suitable change of variables in time, we can write
$$\left\|[(-\Delta)^{\frac{1}{2}}\M(\mathcal{G}_{\varepsilon (t-s)})] \ast\psi\right\|_{L^2_t(L^2_x)}^2\leq \int_{\R}\int_{0}^{+\infty}|\xi|^2\left((1+\xi^2)\frac{\tanh(\xi)}{\xi}\right)e^{-2\varepsilon \eta |\xi|^2\left(1+|\xi|^2\right)}|\widehat{\psi}(\xi)|^2 d\eta d\xi.$$
Now, with the change of variables $\tau=2\varepsilon \eta |\xi|^2\left(1+|\xi|^2\right)$, we have that 
$$\left\|[(-\Delta)^{\frac{1}{2}}\M(\mathcal{G}_{\varepsilon (t-s)})] \ast\psi\right\|_{L^2_t(L^2_x)}^2\leq \int_{\R}\int_{0}^{+\infty}\frac{|\xi|^2\left((1+\xi^2)\frac{\tanh(\xi)}{\xi}\right)}{2\varepsilon |\xi|^2\left(1+\xi^2\right)}e^{-\tau}|\widehat{\psi}(\xi)|^2 d\tau d\xi \le \frac C{\varepsilon} \Vert \psi \Vert_{L^2}^2$$
where we used the fact that the symbol is uniformly bounded, that is, 
\[\frac{|\xi|^2\left((1+\xi^2)\frac{\tanh(\xi)}{\xi}\right)}{2\varepsilon |\xi|^2\left(1+\xi^2\right)}\le \frac C {\varepsilon}.\]
Then, we can deduce the control (\ref{EstimationNoyau2}), from which we find that
$$\left\|\int_{0}^{t}\mathcal{G}_{\varepsilon (t-s)}\ast(-\Delta)^{\frac{1}{2}}\M(u)ds\right\|_{L^2}\leq \frac{C}{\sqrt{\varepsilon}} \underset{\|\psi\|_{L^2}\leq 1}{\sup}\|u\|_{L^2_t(L^2_x)}\|\psi\|_{L^2},$$
and we finally obtain 
\begin{equation}\label{EstimationLineaire}
\left\|\int_{0}^{t}\mathcal{G}_{\varepsilon (t-s)}\ast(-\Delta)^{\frac{1}{2}}\M(u)ds\right\|_{L^\infty_t(L^2_x)}\leq \frac{C}{\sqrt{\varepsilon}} T^{\frac{1}{2}}\|u\|_{L^\infty_t(L^2_x)}.
\end{equation}
\item For the nonlinear term (3) of (\ref{Duhamel}) we can write
\begin{eqnarray*}
\left\|\int_{0}^{t}\mathcal{G}_{\varepsilon (t-s)}\ast\partial_x\left(\frac{u^2}{2}\right)ds\right\|_{L^2}&\leq &\int_{0}^{t}\left\|\mathcal{G}_{\varepsilon (t-s)}\ast\partial_x\left(\frac{u^2}{2}\right)\right\|_{L^2} ds\\
&\leq &C\int_{0}^{t}\|\partial_x\mathcal{G}_{\varepsilon (t-s)}\|_{L^2} \|u^2\|_{L^1} ds.
\end{eqnarray*}
At this point we recall that $\mathcal{G}_{\varepsilon (t-s)}=\mathfrak{g}_{\varepsilon (t-s)}\ast \mathfrak{h}_{\varepsilon (t-s)}$, and we can write
$$\left\|\int_{0}^{t}\mathcal{G}_{\varepsilon (t-s)}\ast\partial_x\left(\frac{u^2}{2}\right)ds\right\|_{L^2}\leq C\int_{0}^{t}\|\mathfrak{g}_{\varepsilon (t-s)}\ast \partial_x\mathfrak{h}_{\varepsilon (t-s)}\|_{L^2} \|u\|_{L^2} \|u\|_{L^2} ds,$$
so we have (recalling that $ \|\partial_x\mathfrak{h}_{\varepsilon (t-s)}\|_{L^2}\leq C(\varepsilon (t-s))^{-\frac{3}{8}}$ and that $\|\mathfrak{g}_{\varepsilon (t-s)}\|_{L^1}\leq C$):
\begin{eqnarray*}
\left\|\int_{0}^{t}\mathcal{G}_{\varepsilon (t-s)}\ast\partial_x\left(\frac{u^2}{2}\right)ds\right\|_{L^2}&\leq &C  \|u\|_{L^\infty_t(L^2_x)} \|u\|_{L^\infty_t(L^2_x)} \int_{0}^{t}\|\mathfrak{g}_{\varepsilon (t-s)}\|_{L^1} \|\partial_x\mathfrak{h}_{\varepsilon (t-s)}\|_{L^2}ds\\
&\leq &C  \|u\|_{L^\infty_t(L^2_x)} \|u\|_{L^\infty_t(L^2_x)} \int_{0}^{t}(\varepsilon (t-s))^{-\frac{3}{8}}ds\\
&\leq &C \varepsilon^{-\frac{3}{8}} t^{\frac{5}{8}}\|u\|_{L^\infty_t(L^2_x)} \|u\|_{L^\infty_t(L^2_x)},
\end{eqnarray*}
and taking the $L^\infty$ norm in the time variable we have 
\begin{equation}\label{EstimationNonLineaire}
\left\|\int_{0}^{t}\mathfrak{g}_{\varepsilon (t-s)}\ast\partial_x\left(\frac{u^2}{2}\right)ds\right\|_{L^\infty_t(L^2_x)}\leq C \varepsilon^{-\frac{3}{8}} T^{\frac{5}{8}}\|u\|_{L^\infty_t(L^2_x)} \|u\|_{L^\infty_t(L^2_x)}.
\end{equation}
\end{itemize}
Now with the estimates (\ref{DatoInicial}), (\ref{EstimationLineaire}) and (\ref{EstimationNonLineaire}) we can easily obtain the existence of a solution $u_\varepsilon$ of the problem (\ref{Duhamel}). Indeed, if we have 
$$\frac{C}{\sqrt{\varepsilon}} T^{\frac{1}{2}}<1, \qquad (C \varepsilon^{-\frac{3}{8}} T^{\frac{5}{8}})\|u_0\|_{L^2}<1\quad\mbox{and}\quad \frac{C}{\sqrt{\varepsilon}} T^{\frac{1}{2}}+(C \varepsilon^{-\frac{3}{8}} T^{\frac{5}{8}})\|u_0\|_{L^2}<1,$$
which is possible for a small time $T$ (recall that $\varepsilon>0$ is fixed), then by the Lemma \ref{Lemma_Existence} we obtain the existence of a mild solution $u_\varepsilon$ of the problem (\ref{Duhamel}) such that $u_\varepsilon \in L^\infty([0,T], L^2(\R))$ for some small time $T>0$. The Proposition \ref{Propo_ExistenceApprox} is now proven. \hfill$\blacksquare$\\

As we can observe, the introduction of the hyperviscosity term $-\varepsilon \mathcal{L}$ allows us to consider the simpler system (\ref{Duhamel}) and the solutions constructed in the above lines satisfy some important properties: 
\begin{Proposition}\label{Propo_Regular}
The solution $u_\varepsilon$ obtained in the Proposition \ref{Propo_ExistenceApprox} is regular. 
\end{Proposition}
\noindent {\bf Proof.} We will first prove here that $u_\varepsilon \in L^\infty([0,T], \dot{H}^{\frac{1}{2}}(\R))$. Indeed, we have by the integral representation formula (\ref{Duhamel}):
\begin{eqnarray}
\|u_\varepsilon\|_{L^\infty_t(\dot{H}^{\frac{1}{2}}_x)}&\leq &\underbrace{\left\|\mathcal{G}_{\varepsilon t}\ast (\varphi_\varepsilon\ast u_0)\right\|_{L^\infty_t(\dot{H}^{\frac{1}{2}}_x)}}_{(1)}+\underbrace{\left\|\int_{0}^{t}\mathcal{G}_{\varepsilon (t-s)}\ast(-\Delta)^{\frac{1}{2}}\M(u_\varepsilon)ds\right\|_{L^\infty_t(\dot{H}^{\frac{1}{2}}_x)}}_{(2)}\label{MildReguliere}\\
&&+\underbrace{\left\|\int_{0}^{t}\mathcal{G}_{\varepsilon (t-s)}\ast\partial_x\left(\frac{u_\varepsilon^2}{2}\right)ds\right\|_{L^\infty_t(\dot{H}^\sigma_x)}}_{(3)},\notag
\end{eqnarray}
and we will study each of these terms separately. 
\begin{itemize}
\item For the initial data (1) in (\ref{MildReguliere}) we have 
$$\left\|\mathcal{G}_{\varepsilon t}\ast (\varphi_\varepsilon\ast u_0)\right\|_{\dot{H}^{\frac{1}{2}}}=\left\|\mathcal{G}_{\varepsilon t}\ast(-\Delta)^{\frac{1}{4}} \varphi_\varepsilon\ast u_0\right\|_{L^2}\leq\|\mathcal{G}_{\varepsilon t}\|_{L^1}\|(-\Delta)^{\frac{1}{4}} \varphi_\varepsilon\|_{L^1}\|u_0\|_{L^2},$$
from which we easily obtain the estimate
$$\left\|\mathcal{G}_{\varepsilon t}\ast (\varphi_\varepsilon\ast u_0)\right\|_{L^\infty_t(\dot{H}_x^\sigma)}\leq C\varepsilon^{-{\frac{1}{2}}}\|u_0\|_{L^2}<+\infty.$$
\item For the linear term (2) in (\ref{MildReguliere}), just as before, by duality we write
\begin{eqnarray*}
\left\|\int_{0}^{t}\mathcal{G}_{\varepsilon (t-s)}\ast(-\Delta)^{\frac{1}{2}}\M(u_\varepsilon)ds\right\|_{\dot{H}^{\frac{1}{2}}}&=&\underset{\|\psi\|_{L^2}\leq 1}{\sup}\left|\int_{0}^{t}\int_{\R} u_\varepsilon[(-\Delta)^{\frac{{\frac{1}{2}}+1}{2}}\M(\mathcal{G}_{\varepsilon (t-s)})] \ast\psi dxds\right|\\
&\leq & \underset{\|\psi\|_{L^2}\leq 1}{\sup}\|u_\varepsilon\|_{L^2_t(L^2_x)}\left\|[(-\Delta)^{\frac{3}{4}}\M(\mathcal{G}_{\varepsilon (t-s)})] \ast\psi\right\|_{L^2_t(L^2_x)}.
\end{eqnarray*}
and we have
$$\left\|[(-\Delta)^{\frac{3}{4}}\M(\mathcal{G}_{\varepsilon (t-s)})] \ast\psi\right\|_{L^2_t(L^2_x)}\leq \frac{C}{\sqrt{\varepsilon}}\|\psi\|_{L^2},$$
since by the Plancherel formula and by the Fubini theorem we write
$$\left\|[(-\Delta)^{\frac{3}{4}}\M(\mathcal{G}_{\varepsilon (t-s)})] \ast\psi\right\|_{L^2_t(L^2_x)}^2\leq \int_{0}^{t}\int_{\R}|\xi|^{3}\left((1+\xi^2)\frac{\tanh(\xi)}{\xi}\right)e^{-2\varepsilon(t-s)|\xi|^2\left(1+|\xi|^2\right)}|\widehat{\psi}(\xi)|^2d\xi ds.$$
In other words, after a change of variables in time, 
$$\left\|[(-\Delta)^{\frac{3}{4}}\M(\mathcal{G}_{\varepsilon (t-s)})] \ast\psi\right\|_{L^2_t(L^2_x)}^2\leq \int_{\R}\int_{0}^{+\infty}|\xi|^{3}\left((1+\xi^2)\frac{\tanh(\xi)}{\xi}\right)e^{-2\varepsilon \eta |\xi|^2\left(1+|\xi|^2\right)}|\widehat{\psi}(\xi)|^2 d\eta d\xi.$$
As we did before, with the change of variables $\tau=2\varepsilon \eta |\xi|^2\left(1+|\xi|^2\right)$, we get 
$$\left\|[(-\Delta)^{\frac{3}{4}}\M(\mathcal{G}_{\varepsilon (t-s)})] \ast\psi\right\|_{L^2_t(L^2_x)}^2\leq \int_{\R}\int_{0}^{+\infty}\frac{|\xi|^{3}\left((1+\xi^2)\frac{\tanh(\xi)}{\xi}\right)}{2\varepsilon |\xi|^2\left(1+\xi^2\right)}e^{-\tau}|\widehat{\psi}(\xi)|^2 d\tau d\xi \leq \frac{C}{\varepsilon}\|\psi\|_{L^2}^2.$$
Thus, we have
$$\left\|\int_{0}^{t}\mathcal{G}_{\varepsilon (t-s)}\ast(-\Delta)^{\frac{1}{2}}\M(u_\varepsilon)ds\right\|_{\dot{H}^\sigma}\leq \frac{C}{\sqrt{\varepsilon}} \underset{\|\psi\|_{L^2}\leq 1}{\sup}\|u_\varepsilon\|_{L^2_t(L^2_x)}\|\psi\|_{L^2},$$
and we finally obtain 
\begin{eqnarray}
\left\|\int_{0}^{t}\mathcal{G}_{\varepsilon (t-s)}\ast(-\Delta)^{\frac{1}{2}}\M(u_\varepsilon)ds\right\|_{L^\infty_t(\dot{H}^{\frac{1}{2}}_x)}&=&\left\|(-\Delta)^{\frac{1}{4}}\int_{0}^{t}\mathcal{G}_{\varepsilon (t-s)}\ast(-\Delta)^{\frac{1}{2}}\M(u_\varepsilon)ds\right\|_{L^\infty_t(L^2_x)}\notag\\
&\leq &\frac{C}{\sqrt{\varepsilon}} T^{\frac{1}{2}}\|u_\varepsilon\|_{L^\infty_t(L^2_x)}<+\infty.\label{EstimationRegulariteLineaire}
\end{eqnarray}
\item Now, for the nonlinear term (3) in (\ref{MildReguliere}) we write
$$\left\|\int_{0}^{t}\mathcal{G}_{\varepsilon (t-s)}\ast\partial_x\left(\frac{u_\varepsilon^2}{2}\right)ds\right\|_{\dot{H}^{\frac{1}{2}}}\leq C\int_{0}^{t}\|(-\Delta)^{\frac{1}{4}}\partial_x\mathcal{G}_{\varepsilon (t-s)}\|_{L^2} \|u_\varepsilon^2\|_{L^1} ds.$$
At this point we recall that $\mathcal{G}_{\varepsilon (t-s)}=\mathfrak{g}_{\varepsilon (t-s)}\ast \mathfrak{h}_{\varepsilon (t-s)}$, and we can write
$$\left\|\int_{0}^{t}\mathcal{G}_{\varepsilon (t-s)}\ast\partial_x\left(\frac{u_\varepsilon^2}{2}\right)ds\right\|_{\dot{H}^{\frac{1}{2}}}\leq C\int_{0}^{t}\|\mathfrak{g}_{\varepsilon (t-s)}\ast (-\Delta)^{\frac{1}{4}}\partial_x\mathfrak{h}_{\varepsilon (t-s)}\|_{L^2} \|u_\varepsilon\|_{L^2} \|u_\varepsilon\|_{L^2} ds,$$
so we have 
$$\left\|\int_{0}^{t}\mathcal{G}_{\varepsilon (t-s)}\ast\partial_x\left(\frac{u_\varepsilon^2}{2}\right)ds\right\|_{\dot{H}^{\frac{1}{2}}}\leq C  \|u_\varepsilon\|_{L^\infty_t(L^2_x)} \|u_\varepsilon\|_{L^\infty_t(L^2_x)} \int_{0}^{t}\|\mathfrak{g}_{\varepsilon (t-s)}\|_{L^1} \|(-\Delta)^{\frac{1}{4}}\partial_x\mathfrak{h}_{\varepsilon (t-s)}\|_{L^2}ds,$$
but since we have the control $\|(-\Delta)^{\frac{1}{4}}\partial_x\mathfrak{h}_{\varepsilon (t-s)}\|_{L^2}\leq C(\varepsilon (t-s))^{-\frac{1}{2}}$ (by a Fourier-based argument), we can write:
\begin{eqnarray*}
\left\|\int_{0}^{t}\mathcal{G}_{\varepsilon (t-s)}\ast\partial_x\left(\frac{u_\varepsilon^2}{2}\right)ds\right\|_{\dot{H}^{\frac{1}{2}}}&\leq &C  \|u_\varepsilon\|_{L^\infty_t(L^2_x)} \|u_\varepsilon\|_{L^\infty_t(L^2_x)} \int_{0}^{t}(\varepsilon (t-s))^{-\frac{1}{2}}ds\\
&\leq &C \varepsilon^{-\frac{1}{2}} t^{\frac{1}{2}}\|u_\varepsilon\|_{L^\infty_t(L^2_x)} \|u_\varepsilon\|_{L^\infty_t(L^2_x)},
\end{eqnarray*}
and taking the $L^\infty$ norm in the time variable we have 
$$\left\|\int_{0}^{t}\mathfrak{g}_{\varepsilon (t-s)}\ast\partial_x\left(\frac{u_\varepsilon^2}{2}\right)ds\right\|_{L^\infty_t(\dot{H}^{\frac{1}{2}}_x)}\leq C \varepsilon^{-\frac{1}{2}} T^{\frac{1}{2}}\|u_\varepsilon\|_{L^\infty_t(L^2_x)} \|u_\varepsilon\|_{L^\infty_t(L^2_x)}<+\infty.$$
\end{itemize}
With all these estimates for the terms in (\ref{MildReguliere}), we have now proven that $u_\varepsilon \in L^{\infty}_t(\dot{H}^\frac{1}{2}_x)$.\\ 

Let us prove now that we have $u_\varepsilon \in L^{\infty}_t(\dot{H}^1_x)$. We write again
\begin{eqnarray*}
\|u_\varepsilon\|_{L^\infty_t(\dot{H}^1_x)}&\leq &\underbrace{\left\|\mathcal{G}_{\varepsilon t}\ast (\varphi_\varepsilon\ast u_0)\right\|_{L^\infty_t(\dot{H}^1_x)}}_{(1)}+\underbrace{\left\|\int_{0}^{t}\mathcal{G}_{\varepsilon (t-s)}\ast(-\Delta)^{\frac{1}{2}}\M(u_\varepsilon)ds\right\|_{L^\infty_t(\dot{H}^1_x)}}_{(2)}\\
&&+\underbrace{\left\|\int_{0}^{t}\mathcal{G}_{\varepsilon (t-s)}\ast\partial_x\left(\frac{u_\varepsilon^2}{2}\right)ds\right\|_{L^\infty_t(\dot{H}^1_x)}}_{(3)},
\end{eqnarray*}
\begin{itemize}
\item Note that the initial data is mollified with the function $\varphi_\varepsilon$, so following the same computations it is easy to prove that 
$$\left\|\mathcal{G}_{\varepsilon t}\ast (\varphi_\varepsilon\ast u_0)\right\|_{L^\infty_t(\dot{H}^1_x)}\leq C\varepsilon^{-1}\|u_0\|_{L^2}<+\infty.$$
\item For the linear term we write, using the commutation properties of the fractional powers of the Laplacian as well as the properties of the operator $\mathcal{M}$:
\begin{eqnarray*}
\left\|\int_{0}^{t}\mathcal{G}_{\varepsilon (t-s)}\ast(-\Delta)^{\frac{1}{2}}\M(u_\varepsilon)ds\right\|_{L^\infty_t(\dot{H}^1_x)}&=&\left\|(-\Delta)^{\frac{1}{2}}\int_{0}^{t}\mathcal{G}_{\varepsilon (t-s)}\ast(-\Delta)^{\frac{1}{2}}\M(u_\varepsilon)ds\right\|_{L^\infty_t(L^2_x)}\\
&=&\left\|(-\Delta)^{\frac{1}{4}}\int_{0}^{t}\mathcal{G}_{\varepsilon (t-s)}\ast(-\Delta)^{\frac{1}{2}}\M((-\Delta)^{\frac{1}{4}}u_\varepsilon)ds\right\|_{L^\infty_t(L^2_x)}.
\end{eqnarray*}
Now, we apply the estimate (\ref{EstimationRegulariteLineaire}) obtained above to write
$$\left\|\int_{0}^{t}\mathcal{G}_{\varepsilon (t-s)}\ast(-\Delta)^{\frac{1}{2}}\M(u_\varepsilon)ds\right\|_{L^\infty_t(\dot{H}^1_x)}\leq \frac{C}{\sqrt{\varepsilon}} T^{\frac{1}{2}}\|(-\Delta)^{\frac{1}{4}}u_\varepsilon\|_{L^\infty_t(L^2_x)}=\frac{C}{\sqrt{\varepsilon}} T^{\frac{1}{2}}\|u_\varepsilon\|_{L^\infty_t(\dot{H}^{\frac{1}{2}}_x)}<+\infty.$$
\item For the nonlinear term, we start writing
$$\left\|\int_{0}^{t}\mathcal{G}_{\varepsilon (t-s)}\ast\partial_x\left(\frac{u_\varepsilon^2}{2}\right)ds\right\|_{\dot{H}^1}=\left\|\int_{0}^{t}(-\Delta)^{\frac{1}{2}}\mathcal{G}_{\varepsilon (t-s)}\ast\partial_x\left(\frac{u_\varepsilon^2}{2}\right)ds\right\|_{L^2},$$
and we have
\begin{eqnarray*}
\left\|\int_{0}^{t}\mathcal{G}_{\varepsilon (t-s)}\ast\partial_x\left(\frac{u_\varepsilon^2}{2}\right)ds\right\|_{\dot{H}^1}&\leq &\int_{0}^{t}\left\|(-\Delta)^{\frac{1}{2}}\mathcal{G}_{\varepsilon (t-s)}\ast\partial_x\left(\frac{u_\varepsilon^2}{2}\right)\right\|_{L^2}ds\\
&\leq &C\int_{0}^{t}\left\|(-\Delta)^{\frac{1}{2}}\partial_x(\mathfrak{g}_{\varepsilon (t-s)}\ast \mathfrak{h}_{\varepsilon (t-s)})\right\|_{L^1}\left\|u_\varepsilon^2\right\|_{L^1}ds.
\end{eqnarray*}
We now write, by the properties of the kernels:
\begin{eqnarray*}
\left\|\int_{0}^{t}\mathcal{G}_{\varepsilon (t-s)}\ast\partial_x\left(\frac{u_\varepsilon^2}{2}\right)ds\right\|_{\dot{H}^1}&\leq &C\int_{0}^{t}\left\|\mathfrak{g}_{\varepsilon (t-s)}\right\|_{L^1} \left\|(-\Delta)^{\frac{1}{2}}\partial_x\mathfrak{h}_{\varepsilon (t-s)}\right\|_{L^2}\left\|u_\varepsilon^2\right\|_{L^1}ds\\
&\leq &C  \|u_\varepsilon\|_{L^\infty_t(L^2_x)} \|u_\varepsilon\|_{L^\infty_t(L^2_x)} \int_{0}^{t}(\varepsilon (t-s))^{-\frac{5}{8}}ds\\
&\leq &C \varepsilon^{-\frac{5}{8}} t^{\frac{3}{8}}\|u_\varepsilon\|_{L^\infty_t(L^2_x)} \|u_\varepsilon\|_{L^\infty_t(L^2_x)},
\end{eqnarray*}
thus taking the $L^\infty$ norm in the space variable we can write
$$\left\|\int_{0}^{t}\mathcal{G}_{\varepsilon (t-s)}\ast\partial_x\left(\frac{u_\varepsilon^2}{2}\right)ds\right\|_{L^\infty_t(\dot{H}^1_x)}\leq C \varepsilon^{-\frac{5}{8}} T^{\frac{3}{8}}\|u_\varepsilon\|_{L^\infty_t(L^2_x)} \|u_\varepsilon\|_{L^\infty_t(L^2_x)}<+\infty.$$
We thus have obtained that $u_\varepsilon\in L^\infty_t(\dot{H}^1_x)$.
\end{itemize}
Now, by repeating the same computations and using the information $u_\varepsilon\in L^\infty_t(\dot{H}^1_x)$ just obtained it is possible to prove that the solution $u_\varepsilon$ is regular, \emph{i.e.} we have $u_\varepsilon\in L^\infty_t(\dot{H}^k_x)$ for some $k>0$ as big as necessary. \hfill$\blacksquare$\\
\begin{Proposition}[\bf Energy inequality]\label{Propo_Energy}
The solutions $u_\varepsilon$ obtained in the previous theorem satisfy the following energy inequality
$$\|u_\varepsilon(t,\cdot)\|_{L^2}^2+\int_{0}^t\|u_\varepsilon(s,\cdot)\|_{\mathcal{N}}^2ds\leq \|\vu_0\|_{L^2}^2,$$
where the functional space $\mathcal{N}(\R)$ is defined by 
$$\|f\|_{\mathcal{N}}^2=\int_{\R}\left|(-\Delta)^{\frac{1}{4}}\M^{\frac{1}{2}}(f)\right|^2dx<+\infty.$$
\end{Proposition}
\noindent {\bf Proof.} Since the solution $u_\varepsilon$ is regular, we can write
$$\int_{\R}(\partial_t u_\varepsilon) u_\varepsilon\, dx = -\varepsilon \int_{\R}\mathcal{L}(u_\varepsilon)u_\varepsilon\,dx-\int_{\R}(-\Delta)^{\frac{1}{2}}\M(u_\varepsilon)u_\varepsilon\,dx+\int_{\R}\partial_x\left(\frac{u_\varepsilon^2}{2}\right)u_\varepsilon\,dx.$$
We remark now that, by an integration by parts,
$$\int_{\R}\partial_x\left(\frac{u_\varepsilon^2}{2}\right)u_\varepsilon\,dx=0$$
We can thus write 
$$\frac{1}{2}\frac{d}{dt}\|u_\varepsilon(t,\cdot)\|_{L^2}^2 + \varepsilon \int_{\R}\left|\mathcal{L}^\frac{1}{2}(u_\varepsilon)\right|^2dx+\int_{\R}\left|(-\Delta)^{\frac{1}{4}}\M^{\frac{1}{2}}(u_\varepsilon)\right|^2dx=0,$$
which implies that 
$$\frac{d}{dt}\|u_\varepsilon(t,\cdot)\|_{L^2}^2 +2\int_{\R}\left|(-\Delta)^{\frac{1}{4}}\M^{\frac{1}{2}}(u_\varepsilon)\right|^2dx\leq 0.$$
Integrating with respect to the time variable, we have 
$$\|u_\varepsilon(t,\cdot)\|_{L^2}^2 +2\int_{0}^t\int_{\R}\left|(-\Delta)^{\frac{1}{4}}\M^{\frac{1}{2}}(u_\varepsilon)\right|^2dxds \leq \|\varphi_\varepsilon\ast u_0\|_{L^2}^2\leq \|u_0\|_{L^2}^2,$$
from which we obtain by definition of the norm $\|\cdot\|_{\mathcal{N}}$:
$$\|u_\varepsilon\|_{L^\infty_t(L^2_x)}^2 +2\int_{0}^t\|u_\varepsilon(s,\cdot)\|_{\mathcal{N}}^2ds \leq  \|u_0\|_{L^2}^2,$$
which is the wished inequality.\hfill$\blacksquare$\\

From this inequality we see that the quantities $L^\infty_t(L^2_x)$ and $L^2_t(\mathcal{N}_x)$ remain uniformly bounded (\emph{i.e.} for any $\varepsilon>0$). Then, using the Duhamel formulation (\ref{Duhamel}) as well as the semi-group properties and the usual arguments for nonlinear PDEs (see \cite{ChNS}, \cite{PGLR2}), we can construct global solutions $u_\varepsilon$ such that $u_\varepsilon \in L^\infty_t(L^2_x) \cap L^2_t(\mathcal{N}_x)$.\\ 

Once he had obtained global solutions $u_\varepsilon \in L^\infty_t(L^2_x) \cap L^2_t(\mathcal{N}_x)$ of the problem (\ref{EqMollified}), we will pass now to the limit $\varepsilon \to 0$. First note that since the sequence $(u_\varepsilon)_{\varepsilon>0}$ is uniformly bounded in the space $L^\infty_t(L^2_x)$, we can apply the Banach-Aloaglu-Bourbaki theorem to consider a subsequence that converges weakly-$*$ to $u\in L^\infty_t(L^2_x)$. Note that we also have the weak limit in $\mathcal{D}'(\R)$ 
$$(-\Delta)^{\frac{1}{2}}\M(u_\varepsilon)\underset{\varepsilon\to 0}{\longrightarrow} (-\Delta)^{\frac{1}{2}}\M(u).$$
However, the convergence of the bilinear term $\partial_x\left(\frac{u_\varepsilon^2}{2}\right)$ requires more details. For this, the following ``interpolation'' lemma will be useful. 
\begin{Lemma}\label{Lemm_InterpolMain}
If $f \in L^\infty_t(L^2_x)\cap L^2_t(\mathcal{N}_x)$ then we have $f\in  L^2_t(\dot{H}_x^{\mathfrak{s}})$ with $0<\mathfrak{s}<\frac{3}{4}$ and 
$$\|f\|_{L^2_t(\dot{H}^\mathfrak{s}_x)}\leq C\|f\|_{L^\infty_t(L^2_x)}^{1-2\mathfrak{s}}\|f\|_{L^2_t(\mathcal{N}_x)}^{2\mathfrak{s}}.$$
\end{Lemma}
\noindent {\bf Proof.} For some parameter $A>0$, to be fixed later, we have 
\begin{eqnarray*}
\|f(t,\cdot)\|_{\dot{H}^\mathfrak{s}}^2&=&\int_{\R}|\xi|^{2\mathfrak{s}}|\widehat{f}(t,\xi)|^2d\xi=\int_{\{|\xi|\leq A\}}|\xi|^{2\mathfrak{s}}|\widehat{f}(t,\xi)|^2d\xi+\int_{\{|\xi|>A\}}|\xi|^{2\mathfrak{s}}|\widehat{f}(t,\xi)|^2d\xi\\
&\leq & A^{2\mathfrak{s}}\int_{\{|\xi|\leq A\}}|\widehat{f}(t,\xi)|^2d\xi+\int_{\{|\xi|>A\}}|\xi|^{2\mathfrak{s}}|\widehat{f}(t,\xi)|^2d\xi. 
\end{eqnarray*}
Noting that $|\xi|^{-\frac{3}{2}}|\xi|\left((1+\xi^2)\frac{\tanh(\xi)}{\xi}\right)^{\frac{1}{2}}>1$ we obtain
$$\|f(t,\cdot)\|_{\dot{H}^\mathfrak{s}}^2\leq  A^{2\mathfrak{s}}\int_{\{|\xi|\leq A\}}|\widehat{f}(t,\xi)|^2d\xi+\int_{\{|\xi|>A\}}|\xi|^{2\mathfrak{s}-\frac{3}{2}}|\xi|\left((1+\xi^2)\frac{\tanh(\xi)}{\xi}\right)^{\frac{1}{2}}|\widehat{f}(t,\xi)|^2d\xi.$$
 Now we write, since $2\mathfrak{s}-\frac{3}{2}<0$, we write
\begin{eqnarray*}
\|f(t,\cdot)\|_{\dot{H}^\mathfrak{s}}^2&\leq & A^{2\mathfrak{s}}\int_{\{|\xi|\leq A\}}|\widehat{f}(t,\xi)|^2d\xi+A^{2\mathfrak{s}-\frac{3}{2}}\int_{\{|\xi|>A\}}|\xi|\left((1+\xi^2)\frac{\tanh(\xi)}{\xi}\right)^{\frac{1}{2}}|\widehat{f}(t,\xi)|^2d\xi\\
&\leq & A^{2\mathfrak{s}}\int_{\R}|\widehat{f}(t,\xi)|^2d\xi+A^{2\mathfrak{s}-\frac{3}{2}}\int_{\R}|\xi|\left((1+\xi^2)\frac{\tanh(\xi)}{\xi}\right)^{\frac{1}{2}}|\widehat{f}(t,\xi)|^2d\xi\\
&\leq & A^{2\mathfrak{s}}\|f(t,\cdot)\|_{L^2}^2+A^{2\mathfrak{s}-\frac{3}{2}}\|f(t,\cdot)\|_{\mathcal{N}}^2.
\end{eqnarray*}
We set now the parameter $A$ such that $A=\frac{\|f(t,\cdot)\|_{\mathcal{N}}^2}{\|f(t,\cdot)\|_{L^2}^2}$ and from the previous inequality we can easily deduce the estimate
$$\|f(t,\cdot)\|_{\dot{H}^\mathfrak{s}}^2\leq C\|f(t,\cdot)\|_{L^2}^{2-4\mathfrak{s}}\|f(t,\cdot)\|_{\mathcal{N}}^{4\mathfrak{s}},$$
thus, taking the $L^2$ norm in the time variable we obtain
$$\|f\|_{L^2_t(\dot{H}^\mathfrak{s}_x)}\leq C\|f\|_{L^\infty_t(L^2_x)}^{1-2\mathfrak{s}}\|f\|_{L^2_t(\mathcal{N}_x)}^{2\mathfrak{s}}<+\infty,$$
which is the wished estimate.\hfill$\blacksquare$\\

With this control at hand, we have that the solutions $u_\varepsilon(t,\cdot)$ belong to the Sobolev space $\dot{H}^{\mathfrak{s}}$ for some $0<\mathfrak{s}<\frac{3}{4}$. In particular, we can now invoque the Rellich-Aubins-Lions theorem that provides us with the compact (local) inclusion $\dot{H}^{\mathfrak{s}}\subset\subset L^q$ for all $1\leq q<\frac{2}{1-2\mathfrak{s}}$ if $0<\mathfrak{s}<\frac{1}{2}$. With this strong (local) inclusion, we can pass to the limit $\varepsilon \to 0$  in the nonlinear term as for the Navier-Stokes equation (see \cite{ChNS}, \cite{PGLR1} or \cite{PGLR2} for the details) and we obtain global weak solutions such that $u \in L^\infty_t(L^2_x)\cap L^2_t(\mathcal{N}_x)$. The Theorem \ref{Theorem_GlobalExistence} is now proven. \hfill$\blacksquare$\\
\begin{Remark}
Since we have $u\in L^\infty_t(L^2_x)\cap L^2_t(\mathcal{N}_x)$, by the Lemma \ref{Lemm_InterpolMain} we have $u\in L^2_t(\dot{H}^{\mathfrak{s}}_x)$ for any $0<\mathfrak{s}<\frac{3}{4}$. 
\end{Remark}
\mysection{Regularity}\label{Secc_Regularity}
In this section, we study the regularity of the weak solutions obtained previously and we will prove Theorem \ref{Theorem_SerrinLinfty}. We will assume here that the initial data $u_0$ satisfies $u_0\in L^2(\R)$ and $u_0\in \dot{H}^{\rho}(\R)$ for some $\rho>0$ big enough and we will assume moreover that $u\in L^\infty_t(L^\infty_x)$. We thus start with the system
\begin{equation}\label{EqRegu}
\begin{cases}
\partial_t u = -(-\Delta)^{\frac{1}{2}}\M(u)+\partial_x\left(\frac{|u|^2}{2}\right),\\[3mm]
u(0,x)=u_0(x), \quad x\in \R,
\end{cases}
\end{equation}
which we rewrite as follows
$$u(t,x)=\mathcal{H}_{t}\ast  u_0(x)+\int_{0}^{t}\mathcal{H}_{(t-s)}\ast\partial_x\left(\frac{u^2}{2}\right)ds,$$
where the kernel $\mathcal{H}_{t}$ is associated to the convolution semi-group $e^{-t(-\Delta)^{\frac{1}{2}}\M}$. \\

We will prove now that we have a first gain of regularity, \emph{i.e.} we will show that  $u\in L^\infty_t(\dot{H}^{\sigma}_x)$ with $0<\sigma=\mathfrak{s}-\frac{1}{4}<\frac{1}{2}$. For this we write
\begin{eqnarray*}
\|u\|_{\dot{H}^\sigma}&\leq &\|\mathcal{H}_{t}\ast  u_0\|_{\dot{H}^\sigma}+\left\|\int_{0}^{t}\mathcal{H}_{(t-s)}\ast\partial_x\left(\frac{u^2}{2}\right)ds\right\|_{\dot{H}^\sigma}\\
&\leq &\|(-\Delta)^{\frac{\sigma}{2}} (\mathcal{H}_{t}\ast u_0)\|_{L^2}+\left\|(-\Delta)^{\frac{\sigma}{2}}\int_{0}^{t}\mathcal{H}_{(t-s)}\ast\partial_x\left(\frac{u^2}{2}\right)ds\right\|_{L^2}.
\end{eqnarray*}
The first term above is easy to handle as we have
\begin{equation}\label{DatoInicialRegu}
\|(-\Delta)^{\frac{\sigma}{2}} (\mathcal{H}_{t}\ast u_0)\|_{L^2}\leq \|\mathcal{H}_{t}\|_{L^1}\|(-\Delta)^{\frac{\sigma}{2}}  u_0\|_{L^2}\leq C\|u_0\|_{\dot{H}^\sigma}<+\infty.
\end{equation}
For the bilinear term we write
$$\left\|(-\Delta)^{\frac{\sigma}{2}}\int_{0}^{t}\mathcal{H}_{(t-s)}\ast\partial_x\left(\frac{u^2}{2}\right)ds\right\|_{L^2}\leq C\left\|\int_{0}^{t}(-\Delta)^{\frac{\sigma-\mathfrak{s}}{2}}\partial_x\mathcal{H}_{(t-s)}\ast (-\Delta)^{\frac{\mathfrak{s}}{2}}(u^2)ds\right\|_{L^2},$$
and by duality we have 
$$\left\|(-\Delta)^{\frac{\sigma}{2}}\int_{0}^{t}\mathcal{H}_{(t-s)}\ast\partial_x\left(\frac{u^2}{2}\right)ds\right\|_{L^2}\leq C\underset{\|\psi\|_{L^2}\leq 1}{\sup}\left|\int_{\R} \left(\int_{0}^{t}(-\Delta)^{\frac{\sigma-\mathfrak{s}}{2}}\partial_x\mathcal{H}_{(t-s)}\ast (-\Delta)^{\frac{\mathfrak{s}}{2}}(u^2)ds\right)\psi dx\right|.$$
By extending the kernel $\mathcal{H}_{(t-s)}$ by $0$ if $s>t$, we can write
$$\left\|(-\Delta)^{\frac{\sigma}{2}}\int_{0}^{t}\mathcal{H}_{(t-s)}\ast\partial_x\left(\frac{u^2}{2}\right)ds\right\|_{L^2}\leq C\underset{\|\psi\|_{L^2}\leq 1}{\sup}\left|\int_{\R} \left(\int_{0}^{+\infty}(-\Delta)^{\frac{\sigma-\mathfrak{s}}{2}}\partial_x\mathcal{H}_{(t-s)}\ast (-\Delta)^{\frac{\mathfrak{s}}{2}}(u^2)ds\right)\psi dx\right|,$$
Then, by the properties of the convolution, we obtain
$$\left\|(-\Delta)^{\frac{\sigma}{2}}\int_{0}^{t}\mathcal{H}_{(t-s)}\ast\partial_x\left(\frac{u^2}{2}\right)ds\right\|_{L^2}\leq C\underset{\|\psi\|_{L^2}\leq 1}{\sup}\left|\int_{0}^{+\infty}\int_{\R}  (-\Delta)^{\frac{\mathfrak{s}}{2}}(u^2)[(-\Delta)^{\frac{\sigma-\mathfrak{s}}{2}}\partial_x\mathcal{H}_{(t-s)}\ast\psi] dxds\right|,$$
Finally, from the Cauchy-Schwarz inequality, we have that 
\begin{eqnarray}
\left\|(-\Delta)^{\frac{\sigma}{2}}\int_{0}^{t}\mathcal{H}_{(t-s)}\ast\partial_x\left(\frac{u^2}{2}\right)ds\right\|_{L^2} C \left\|u^2\right\|_{L^2_t(\dot{H}^{\mathfrak{s}}_x)} \;\underset{\|\psi\|_{L^2}\leq 1}{\sup}\left\|(-\Delta)^{\frac{\sigma-\mathfrak{s}}{2}}\partial_x\mathcal{H}_{(t-s)}\ast\psi\right\|_{L^2_t(L^2_x)}.\label{GainRegulariteSolFaible}
\end{eqnarray}
We will study these two terms separately. For the first one we have, by the Kato-Ponce inequalities (\ref{KatoPonce_Ineq}):
$$\left\|u^2(t,\cdot)\right\|_{\dot{H}^{\mathfrak{s}}}\leq C\|u(t,\cdot)\|_{\dot{H}^\mathfrak{s}}\|u(t,\cdot)\|_{L^\infty}.$$
Now, taking the $L^2$ norm in the time variable, we obtain (recall that we assumed $u\in L^\infty_t(L^\infty_x)$)
\begin{equation}\label{KatoPonce}
\left\|u^2\right\|_{L^2_t(\dot{H}^{\mathfrak{s}}_x)}\leq C\|u\|_{L^2_t(\dot{H}^\mathfrak{s}_x)}\|u\|_{L^\infty_t(L^\infty_x)}.
\end{equation}
For the second term of (\ref{GainRegulariteSolFaible}), we write
$$\left\|[(-\Delta)^{\frac{\sigma-\mathfrak{s}}{2}}\partial_x\mathcal{H}_{(t-s)})] \ast\psi\right\|_{L^2_t(L^2_x)}^2\leq \int_{0}^{+\infty}\int_{\R}|\xi|^{2+2(\sigma-\mathfrak{s})}e^{-2(t-s)|\xi|\sqrt{(1+\xi^2)\frac{\tanh(\xi)}{\xi}}}|\widehat{\psi}(\xi)|^2d\xi ds.$$
Since $\sigma=\mathfrak{s}-\frac{1}{4}>0$, we obtain
$$\left\|[(-\Delta)^{\frac{\sigma-\mathfrak{s}}{2}}\partial_x\mathcal{H}_{(t-s)})] \ast\psi\right\|_{L^2_t(L^2_x)}^2 \leq \int_{0}^{+\infty}\int_{\R}|\xi|^{\frac{3}{2}}e^{-2(t-s)|\xi|\sqrt{(1+\xi^2)\frac{\tanh(\xi)}{\xi}}}|\widehat{\psi}(\xi)|^2d\xi ds.$$
Setting $\eta=2(t-s)|\xi|\sqrt{(1+\xi^2)\frac{\tanh(\xi)}{\xi}}$, by a change of variables we have
$$\left\|[(-\Delta)^{\frac{\sigma-\mathfrak{s}}{2}}\partial_x\mathcal{H}_{(t-s)})] \ast\psi\right\|_{L^2_t(L^2_x)}^2\leq\int_{0}^{+\infty}\int_{\R}\frac{|\xi|^{\frac{3}{2}}}{2|\xi|\sqrt{(1+\xi^2)\frac{\tanh(\xi)}{\xi}}}e^{-\eta}|\widehat{\psi}(\xi)|^2d\xi d\eta.$$
We remark again that the symbol $\frac{|\xi|^{\frac{3}{2}}}{2|\xi|\sqrt{(1+\xi^2)\frac{\tanh(\xi)}{\xi}}}$ is uniformly bounded, so we can write 
\begin{equation}\label{Serrin}
\left\|[(-\Delta)^{\frac{\sigma-\mathfrak{s}}{2}}\partial_x\mathcal{H}_{(t-s)})] \ast\psi\right\|_{L^2_t(L^2_x)}^2\leq C\int_{\R}\int_{0}^{+\infty}e^{-\eta} d\eta\, |\widehat{\psi}(\xi)|^2d\xi\leq C\|\psi\|_{L^2}^2<+\infty,
\end{equation}
since we assumed $\|\psi\|_{L^2}\leq 1$.\\

With the estimates (\ref{KatoPonce}) and (\ref{Serrin}) at hand, we come back to (\ref{GainRegulariteSolFaible}), to obtain
\begin{equation}\label{BilinealRegu}
\left\|(-\Delta)^{\frac{\sigma}{2}}\int_{0}^{t}\mathcal{H}_{(t-s)}\ast\partial_x\left(\frac{u^2}{2}\right)ds\right\|_{L^2}\leq  C \|u\|_{L^2_t(\dot{H}^\mathfrak{s}_x)}\|u\|_{L^\infty_t(L^\infty_x)}<+\infty.
\end{equation}

Now, with the controls (\ref{DatoInicialRegu}) and (\ref{BilinealRegu}) we have proven that the weak solution $u$ of the system (\ref{EqRegu}) belongs to the space $L^\infty_t(\dot{H}^\sigma_x)$ with $0<\sigma<\frac{1}{2}$.\\

With the previous step, where we obtained under the hypothesis $u\in L^\infty_t(L^\infty_x)$ a first gain of regularity $u\in L^\infty_t(\dot{H}^\sigma_x)$ with $0<\sigma<\frac{1}{2}$, we will now prove a second gain of regularity: indeed we will obtain that $u\in L^\infty_t(\dot{H}^\alpha_x)$ with $0<\alpha<\frac{7}{4}$.\\

Note that the initial data is regular by hypothesis, so the same computations performed in (\ref{DatoInicialRegu}) apply here. We will focus now in the nonlinear term and we have 
$$\left\|(-\Delta)^{\frac{\alpha}{2}}\int_{0}^{t}\mathcal{H}_{(t-s)}\ast\partial_x\left(\frac{u^2}{2}\right)ds\right\|_{L^2}\leq C\left\|\int_{0}^{t}(-\Delta)^{\frac{\alpha-(\mathfrak{s}+\sigma-\frac{1}{2})}{2}}\partial_x\mathcal{H}_{(t-s)}\ast (-\Delta)^{\frac{\mathfrak{s}+\sigma-\frac{1}{2}}{2}}(u^2)ds\right\|_{L^2},$$
recall that $\mathfrak{s}<\frac{3}{4}$ and $\sigma<\frac{1}{2}$, so we have $\mathfrak{s}+\sigma<\frac{5}{2}$.\\

Following the same ideas displayed to obtain (\ref{GainRegulariteSolFaible}), we can write
\begin{eqnarray}
\left\|(-\Delta)^{\frac{\alpha}{2}}\int_{0}^{t}\mathcal{H}_{(t-s)}\ast\partial_x\left(\frac{u^2}{2}\right)ds\right\|_{L^2}&\leq & C \left\|u^2\right\|_{L^2_t(\dot{H}^{\mathfrak{s}+\sigma-\frac{1}{2}}_x)} \notag\\
&&\;\underset{\|\psi\|_{L^2}\leq 1}{\sup}\left\|(-\Delta)^{\frac{\alpha-(\mathfrak{s}+\sigma-\frac{1}{2})}{2}}\partial_x\mathcal{H}_{(t-s)}\ast\psi\right\|_{L^2_t(L^2_x)}.\label{GainRegulLoiProduit}
\end{eqnarray}
Now, for the first term above we use the usual product laws in Sobolev spaces, see the estimate (\ref{LawProduct_Ineq}) above, to write
$$\left\|u^2\right\|_{\dot{H}^{\mathfrak{s}+\sigma-\frac{1}{2}}}\leq C\|u\|_{\dot{H}^{\mathfrak{s}}}\|u\|_{\dot{H}^{\sigma}},$$
and taking the $L^2$ norm in the time variable one gets
$$\left\|u^2\right\|_{L^2_t(\dot{H}^{\mathfrak{s}+\sigma-\frac{1}{2}}_x)}\leq C\|u\|_{L^2_t(\dot{H}^{\mathfrak{s}}_x)}\|u\|_{L^\infty_t(\dot{H}^{\sigma}_x)}.$$
For the second term of (\ref{GainRegulLoiProduit}), we proceed as before as long as $\alpha-(\mathfrak{s}+\sigma-\frac{1}{2})=\frac{1}{4}$, which lead us to 
$$\left\|(-\Delta)^{\frac{\alpha-(\mathfrak{s}+\sigma-\frac{1}{2})}{2}}\partial_x\mathcal{H}_{(t-s)}\ast\psi\right\|_{L^2_t(L^2_x)}^2\leq \int_{0}^{+\infty}\int_{\R}|\xi|^{\frac{3}{2}}e^{-2(t-s)|\xi|\sqrt{(1+\xi^2)\frac{\tanh(\xi)}{\xi}}}|\widehat{\psi}(\xi)|^2d\xi ds,$$
and following the same ideas that lead us to (\ref{Serrin}) we obtain 
$$\left\|(-\Delta)^{\frac{\alpha-(\mathfrak{s}+\sigma-\frac{1}{2})}{2}}\partial_x\mathcal{H}_{(t-s)}\ast\psi\right\|_{L^2_t(L^2_x)}^2\leq C\|\psi\|_{L^2}^2<+\infty.$$
Note now that we imposed the condition $\alpha=\mathfrak{s}+\sigma-\frac{1}{4}$, but since $\mathfrak{s}<\frac{3}{4}$ and $\sigma<\frac{1}{2}$ we have $\alpha<1$, we obtain now that 
$u\in L^\infty_t(\dot{H}^\alpha_x)$ with $\alpha<1$.\\

By iteration, we deduce that $u\in L^\infty_t(\dot{H}^\rho_x)$ where the maximal value of the parameter $\rho$ is fixed by the information available over the initial data $u_0$. Theorem \ref{Theorem_SerrinLinfty} is now proven. \hfill $\blacksquare$\\

\begin{Remark}\label{Rem_Linfty}
Note that it is possible to go one step further in the regularity study, as from the information $u\in L^\infty_t(\dot{H}^\rho_x)$ for $\rho$ big enough, we can deduce by the Sobolev embeddings the information $\partial_xu\in L^\infty_t(L^\infty_x)$.
\end{Remark}
\mysection{Uniqueness}\label{Secc_Uniqueness}
One interesting consequence of regularity of weak solutions is the uniqueness and we will prove now Theorem \ref{Theorem_Uniqueness}. Recall that by hypothesis we have $u\in L^\infty_t(L^\infty_x)$ thus by Theorem \ref{Theorem_SerrinLinfty} we also have $u\in L^\infty_t(\dot{H}^\rho_x)$ where $\rho>0$ is given by the information over the initial data. \\

Assume now that $v$ is also a weak solution of the equation (\ref{EqIntro}) arising from the same initial data $u_0$ and such that $v\in L^\infty_t(L^\infty_x)$ (we also have  $v\in L^\infty_t(\dot{H}^\rho_x)$).\\

We consider the variable $w=u-v$, since $u$ and $v$ have the same initial data $u_0$ we have 
$$w(0,x)=u(0,x)-v(0,x)=u_0(x)-u_0(x)=0.$$ 
If we study the dynamics of $w$ we can write 
$$\partial_t w=\partial_t u-\partial_t v=\left(-(-\Delta)^{\frac{1}{2}}\mathcal{M}(u)+\partial_x(\frac{u^2}{2})\right)-\left(-(-\Delta)^{\frac{1}{2}}\mathcal{M}(v)+\partial_x(\frac{v^2}{2})\right),$$
and by the linearity of the operators we obtain
\begin{eqnarray*}
\partial_t w&=&-(-\Delta)^{\frac{1}{2}}\mathcal{M}(u-v)+\partial_x(\frac{u^2-v^2}{2})\\
&=&-(-\Delta)^{\frac{1}{2}}\mathcal{M}(w)+\partial_x(\frac{w(u^2+v^2)}{2}).
\end{eqnarray*}
Since we have enough regularity for the variables $u,v,w$, we have the expression
$$\int_{\R}[\partial_t w]w\,dx=-\int_{R}[(-\Delta)^{\frac{1}{2}}\mathcal{M}(w)]w\,dx+\int_{\R}[\partial_x(\frac{w(u^2+v^2)}{2})]w\,dx,$$
which can we rewritten as follows
\begin{eqnarray}
\frac{1}{2}\frac{d}{dt}\|w(t,\cdot)\|_{L^2}^2+\|(-\Delta)^{\frac{1}{4}}\mathcal{M}^\frac{1}{2}(w)\|_{L^2}^2&=&\int_{\R}[\partial_x(\frac{w(u^2+v^2)}{2})]w\,dx\notag\\
\frac{1}{2}\frac{d}{dt}\|w(t,\cdot)\|_{L^2}^2&\leq&\int_{\R}[\partial_x(\frac{w(u^2+v^2)}{2})]w\,dx\notag\\
&\leq &\frac{1}{2}\int_{\R}[(\partial_x w)(u^2+v^2)]w\,dx+\frac{1}{2}\int_{\R}[\partial_x(u^2+v^2)]w^2\,dx.\quad \label{Gron1}
\end{eqnarray}
In particular, by an integration by parts we get
$$\int_{\R}[(\partial_x w)(u^2+v^2)]w\,dx=-\int_{\R}[(\partial_x w)(u^2+v^2)]w\,dx,$$
from which we deduce that 
$$\int_{\R}[(\partial_x w)(u^2+v^2)]w\,dx=0,$$
and thus we can rewrite (\ref{Gron1}) as
$$\frac{1}{2}\frac{d}{dt}\|w(t,\cdot)\|_{L^2}^2\leq\frac{1}{2}\int_{\R}[\partial_x(u^2+v^2)]w^2\,dx.$$
Since $u,v\in L^\infty_t(L^\infty_x)$ and $\partial_xu, \partial_xv\in L^\infty_t(L^\infty_x)$ (recall Remark \ref{Rem_Linfty}) we obtain
$$\frac{1}{2}\frac{d}{dt}\|w(t,\cdot)\|_{L^2}^2\leq \left(\|u\|_{L^\infty_t(L^\infty_x)}\|\partial_xu\|_{L^\infty_t(L^\infty_x)}+\|v\|_{L^\infty_t(L^\infty_x)}\|\partial_xv\|_{L^\infty_t(L^\infty_x)}\right)\|w(t,\cdot)\|_{L^2}^2,$$
now, integrating in the time variable we obtain
$$\|w(t,\cdot)\|_{L^2}^2\leq C\left(\|u\|_{L^\infty_t(L^\infty_x)}\|\partial_xu\|_{L^\infty_t(L^\infty_x)}+\|v\|_{L^\infty_t(L^\infty_x)}\|\partial_xv\|_{L^\infty_t(L^\infty_x)}\right)\int_{0}^t\|w(s,\cdot)\|_{L^2}^2ds+\|w(0,\cdot)\|_{L^2}^2,$$
but since $w(0,\cdot)=0$, we have the estimate
$$\|w(t,\cdot)\|_{L^2}^2\leq C\left(\|u\|_{L^\infty_t(L^\infty_x)}\|\partial_xu\|_{L^\infty_t(L^\infty_x)}+\|v\|_{L^\infty_t(L^\infty_x)}\|\partial_xv\|_{L^\infty_t(L^\infty_x)}\right)\int_{0}^t\|w(s,\cdot)\|_{L^2}^2ds,$$
from which we deduce, by the Gr\"onwall lemma, that $w(t,\cdot)=0$. We have thus obtained the wished uniqueness result and the Theorem \ref{Theorem_Uniqueness} is now proven. \hfill $\blacksquare$

\mysection{Local in time results}\label{Secc_Local}
Note that the previous results are global in time, however if we work in a bounded time interval, say $[0,T]$ for some fixed $0<T<+\infty$, then the hypotheses $u\in L^\infty_tL^\infty_x$ (used for regularity and uniqueness) can be relaxed. Indeed, we first remark that the Lemma \ref{Lemm_InterpolMain} admits the following modification:
\begin{Lemma}\label{Lema34}
If $f \in L^\infty_t(L^2_x)\cap L^2_t(\mathcal{N}_x)$ then we have $f\in  L^2([0,T], \dot{H}^{\frac{3}{4}}(\mathbb{R}))$.
\end{Lemma}
\noindent {\bf Proof.} Following the same ideas as in the proof of Lemma \ref{Lemm_InterpolMain}, we write
\begin{eqnarray*}
\|f(t,\cdot)\|_{\dot{H}^{\frac{3}{4}}}^2&=&\int_{\{|\xi|\leq 1\}}|\xi|^{\frac{3}{2}}|\widehat{f}(t,\xi)|^2d\xi+\int_{\{|\xi|>1\}}|\xi|^{\frac{3}{2}}|\widehat{f}(t,\xi)|^2d\xi\\
&\leq & \int_{\{|\xi|\leq 1\}}|\widehat{f}(t,\xi)|^2d\xi+\int_{\{|\xi|>1\}}|\xi|^{\frac{3}{2}}|\widehat{f}(t,\xi)|^2d\xi,
\end{eqnarray*}
and since $|\xi|^{-\frac{3}{2}}|\xi|\left((1+\xi^2)\frac{\tanh(\xi)}{\xi}\right)^{\frac{1}{2}}>1$ if $|\xi|>1$ we obtain
$$\|f(t,\cdot)\|_{\dot{H}^{\frac{3}{4}}}^2\leq \int_{\{|\xi|\leq 1\}}|\widehat{f}(t,\xi)|^2d\xi+\int_{\{|\xi|>1\}}|\xi|\left((1+\xi^2)\frac{\tanh(\xi)}{\xi}\right)^{\frac{1}{2}}|\widehat{f}(t,\xi)|^2d\xi.$$
Now we write
\begin{eqnarray*}
\|f(t,\cdot)\|_{\dot{H}^{\frac{3}{4}}}^2&\leq &\int_{\R}|\widehat{f}(t,\xi)|^2d\xi+\int_{\R}|\xi|\left((1+\xi^2)\frac{\tanh(\xi)}{\xi}\right)^{\frac{1}{2}}|\widehat{f}(t,\xi)|^2d\xi\\
&\leq & \|f(t,\cdot)\|_{L^2}^2+\|f(t,\cdot)\|_{\mathcal{N}}^2.
\end{eqnarray*}
Taking the $L^2$ norm in the time variable we obtain
$$\|f\|_{L^2_t\left(\dot{H}^{\frac{3}{4}}_x\right)}\leq T^{\frac{1}{2}}\|f\|_{L^\infty_t\left(L^2_x\right)}+\|f\|_{L^2_t(\mathcal{N}_x)},$$
which is a bounded quantity since $T<+\infty$.\hfill$\blacksquare$\\

This result is interesting as it is an improvement to the endpoint $\frac{3}{4}$ of the Lemma \ref{Lemm_InterpolMain} but it is not enough to our purposes. We will now see that under a mild hypothesis,  any weak solution $u\in L^\infty_t(L^2_x)\cap L^2_t(\mathcal{N}_x)$ obtained in the Theorem \ref{Theorem_GlobalExistence} is bounded in an interval $[0,T]$:
\begin{Proposition}
Let $u_0\in L^2(\R)\cap L^\infty(\R)$ be an initial data and consider the associated weak solution $u\in L^\infty_t(L^2_x)\cap L^2_t(\mathcal{N}_x)$ to the equation (\ref{EqIntro}) obtained via the Theorem \ref{Theorem_GlobalExistence}. For a bounded time interval $[0,T]$ with $0<T<+\infty$, if we assume that $u\in L^4([0,T], H^{\frac{3}{4}+\epsilon}(\mathbb{R}))$ for some $0<\epsilon\ll 1$, then we have $u\in L^\infty([0,T], L^\infty(\mathbb{R}))$.
\end{Proposition}

\begin{Remark} If we work in a bounded in time interval, then we have $u\in L^2_t\left(L^2_x\right)$ and with the information $u\in L^2_t\left(\dot{H}^{\frac{3}{4}}_x\right)$ given by the Lemma \ref{Lema34} we obtain that $u\in L^2_t(H^{\frac{3}{4}}_x)$. Note also that we do have $u\in L^2_t\left(H^{\frac{3}{4}}_x\right)$ thus the hypothesis $u\in L^4_t\left(H^{\frac{3}{4}+\epsilon}_x\right)$ for some $0<\epsilon\ll 1$ requires an additional integrability condition in the time variable and some extra regularity information in the space variable, and these conditions are slightly weaker than the condition $u\in L^\infty_t(L^\infty_x)$ asked in the Theorem \ref{Theorem_SerrinLinfty} and Theorem \ref{Theorem_Uniqueness}.

\end{Remark}
\noindent {\bf Proof.} We consider the integral formula of the equation (\ref{EqIntro})
$$u(t,x)=\mathcal{H}_{t}\ast  u_0(x)+\int_{0}^{t}\mathcal{H}_{(t-s)}\ast\partial_x\left(\frac{u^2}{2}\right)ds,$$
and we write
$$\|u(t,\cdot)\|_{L^\infty}\leq \|\mathcal{H}_{t}\ast  u_0\|_{L^\infty}+\int_{0}^{t}\left\|\mathcal{H}_{(t-s)}\ast\partial_x\left(\frac{u^2}{2}\right)\right\|_{L^\infty}ds,$$
thus, by the Young inequalities for the convolution we obtain
\begin{eqnarray*}
\|u(t,\cdot)\|_{L^\infty}&\leq &\|\mathcal{H}_{t}\|_{L^1} \|u_0\|_{L^\infty}+C\int_{0}^{t}\left\|\frac{\partial_x}{(-\Delta)^{\frac{3}{8}+\frac{\epsilon}{2}}}\mathcal{H}_{(t-s)}\right\|_{L^2}\|(-\Delta)^{\frac{3}{8}+\frac{\epsilon}{2}}(u^2)\|_{L^2}ds\\
&\leq &  \|\mathcal{H}_{t}\|_{L^1} \|u_0\|_{L^\infty}+C\int_{0}^{t}\left\|(-\Delta)^{\frac{1}{8}-\frac{\epsilon}{2}}\mathcal{H}_{(t-s)}\right\|_{L^2}\|u^2\|_{H^{\frac{3}{4}+\epsilon}}ds,
\end{eqnarray*}
and by the Cauchy-Schwarz inequality in the time variable we have
\begin{eqnarray}
\|u(t,\cdot)\|_{L^\infty}&\leq &\|\mathcal{H}_{t}\|_{L^1} \|u_0\|_{L^\infty}+\|(-\Delta)^{\frac{1}{8}-\frac{\epsilon}{2}}\mathcal{H}_{(t-s)}\|_{L^2_t(L^2_x)}\|u^2\|_{L^2_t(H^{\frac{3}{4}+\epsilon}_x)}\notag\\
&\leq &C\|u_0\|_{L^\infty}+\|(-\Delta)^{\frac{1}{8}-\frac{\epsilon}{2}}\mathcal{H}_{(t-s)}\|_{L^2_t(L^2_x)}\|u^2\|_{L^2_t(H^{\frac{3}{4}+\epsilon}_x)},\label{LInfty0}
\end{eqnarray}
since $\|\mathcal{H}_{t}\|_{L^1}\leq C$ (recall that $\widehat{\mathcal{H}_{t}}(\xi)=e^{-t|\xi|\mathfrak{m}(\xi)}$).\\

We remark now that $\|(-\Delta)^{\frac{1}{8}-\frac{\epsilon}{2}}\mathcal{H}_{(t-s)}\|_{L^2_t(L^2_x)}<+\infty$ since we have
$$\|(-\Delta)^{\frac{1}{8}-\frac{\epsilon}{2}}\mathcal{H}_{(t-s)}\|_{L^2_t(L^2_x)}^2=\int_{0}^t\int_{\R}|\xi|^{\frac{1}{2}-2\epsilon}e^{-2(t-s)|\xi|\mathfrak{m}(\xi)}d\xi ds,$$
by the Fubini theorem it follows
$$\|(-\Delta)^{\frac{1}{8}-\frac{\epsilon}{2}}\mathcal{H}_{(t-s)}\|_{L^2_t(L^2_x)}^2\leq \int_{\R}\int_{0}^{+\infty}|\xi|^{\frac{1}{2}-2\epsilon}e^{-2(t-s)|\xi|\mathfrak{m}(\xi)} dsd\xi,$$
and setting $\tau=2(t-s)|\xi|\mathfrak{m}(\xi)$, by a change of variables we obtain
\begin{eqnarray*}
\|(-\Delta)^{\frac{1}{8}-\frac{\epsilon}{2}}\mathcal{H}_{(t-s)}\|_{L^2_t(L^2_x)}^2&\leq &C\int_{\R}\int_{0}^{+\infty}\frac{|\xi|^{\frac{1}{2}-2\epsilon}}{|\xi|\mathfrak{m}(\xi)}e^{-\tau} d \tau d\xi\leq C\int_{\R}\frac{|\xi|^{-\frac{1}{2}-2\epsilon}}{\mathfrak{m}(\xi)}d\xi\\
&\leq&C\int_{\{|\xi|\leq 1\}}\frac{|\xi|^{-\frac{1}{2}-2\epsilon}}{\mathfrak{m}(\xi)}d\xi+C\int_{\{|\xi|>1\}}\frac{|\xi|^{-\frac{1}{2}-2\epsilon}}{\mathfrak{m}(\xi)}d\xi,
\end{eqnarray*}
recalling that $\mathfrak{m}(\xi)\sim C$ if $|\xi|\leq 1$ and $\mathfrak{m}(\xi)\sim |\xi|^{\frac{1}{2}}$ if $|\xi|>1$, we can write
\begin{eqnarray*}
\|(-\Delta)^{\frac{1}{8}-\frac{\epsilon}{2}}\mathcal{H}_{(t-s)}\|_{L^2_t(L^2_x)}^2&\leq&C\int_{\{|\xi|\leq 1\}}|\xi|^{-\frac{1}{2}-2\epsilon}d\xi+C\int_{\{|\xi|>1\}}|\xi|^{-1-2\epsilon}d\xi<+\infty,
\end{eqnarray*}
as long as $0<\epsilon\ll 1$ is small enough.\\

Thus, coming back to (\ref{LInfty0}), we obtain 
$$\|u(t,\cdot)\|_{L^\infty}\leq C \|u_0\|_{L^\infty}+C\|u^2\|_{L^2_t(H^{\frac{3}{4}+\epsilon}_x)},$$
but since by the Kato-Ponce inequalities (\ref{KatoPonce_Ineq}) we have
\begin{eqnarray*}
\|u^2\|_{L^2_t(H^{\frac{3}{4}+\epsilon}_x)}&=&\left\|\|u^2\|_{H^{\frac{3}{4}+\epsilon}_x}\right\|_{L^2_t}\leq C\left\|\|u\|_{H^{\frac{3}{4}+\epsilon}_x}\|u\|_{L^\infty_x}\right\|_{L^2_t}\\
&\leq &C\left\|\|u\|_{H^{\frac{3}{4}+\epsilon}_x}\|u\|_{H^{\frac{3}{4}+\epsilon}_x}\right\|_{L^2_t}\leq C\|u\|_{L^4_t(H^{\frac{3}{4}+\epsilon}_x)}^2<+\infty,
\end{eqnarray*}
where in the last line we used the Sobolev embedding $H^{\frac{3}{4}+\epsilon}(\R)\subset L^\infty(\R)$  and the Cauchy-Schwarz inequality in the time variable.\\

Thus, we have finally obtained that 
$$\|u(t,\cdot)\|_{L^\infty}\leq C \|u_0\|_{L^\infty}+C\|u\|_{L^4_t(H^{\frac{3}{4}+\epsilon}_x)}^2<+\infty,$$
for all $0<t<T$ and we deduce that $u\in L^\infty([0,T], L^\infty(\mathbb{R}))$. \hfill$\blacksquare$\\

From this information, we can apply Theorems \ref{Theorem_SerrinLinfty} and \ref{Theorem_Uniqueness} to obtain a gain of regularity and to deduce uniqueness results. \\

\noindent {\bf Acknowledgment.} The work of the second author has been partially supported by the project CRISIS (ANR-20-CE40-0020-01), operated by the French National Research Agency (ANR) and by the project HiCE, operated by the Université Claude Bernard Lyon 1. The authors thank the support of the program MATH-AMSUD 23-MATH-18.

\end{document}